\documentclass{amsart}
\usepackage{graphicx}
\theoremstyle{definition}

\theoremstyle{remark}

\begin{document}

\title{Influence of singularities in the problem on the
Convergence Exponent of Multidimensional Terry's problem}

\author{I.Sh.Jabbarov}
\address{AZ2000 Haydar Aliyev avenue, 159, Ganja State University, Azerbaijan}
\email{ilgar\_j@rambler.ru}


\begin{abstract}
In the paper it is considered the question on Convergence Exponent
of the special integral of multidimensional Terry's problem. The
problem considered in the article consisted in investigation of
singularities influence to the value of convergence exponent
of the special integral.

\end{abstract}
\maketitle

\begin{center}
\noindent \textbf{1. Introduction}
\end{center}

In Multidimensional Analysis many questions of the theory of trigonometrical integrals lead to the investigation of singularities of some mappings. They have crucial importance when questions on asymptotics of oscillatory integrals are studying (see [1]). Influence of singularities can vary depending on character of the problem connected with trigonometrical integrals. In the simplest form singularities arise (see [2]) when one consider the question on the estimation of the integral
\[\int _{a}^{b}e^{2i\pi f(x)} dx \]
where $a,b\in \mathbb{R}$, $f(x)$ is a differentiable function. By using of high derivatives one can estimate the integral in neighborhood of the singular points (i. e the points where $f^{k}(x)=0, k\ge1$ ). In the works [3-6], [8] multiple trigonometrical integrals were studied when the authors generalized Vinogradov's method to the multi-dimensional case. Firstly Chubaricov V. N. gave an estimation for multiple trigonometrical integrals (see [5-6]).

The question on the convergence exponent of Terry's problem in one-dimensional case was considered in [7]. Multidimensional case was studied in [4, 8]. Let the polynomial $F(\bar{x})$ be defined by an equality

\[F(\bar{x})=\sum_{j=1}^{N}\alpha_{j}\gamma_{j}(\bar{x}); \bar{x}=(x_{1},x_{2},...,x_{r}),\]
where $\gamma_{j}(\bar{x})= x_{1}^{k_{j1}} x_{2}^{k_{j2}}\cdots x_{r}^{k_{jr}}$ are monomials of degree $k(j)=k_{j1}+ k_{j2}+\cdots +k_{jr}$, moreover, the polynomial has not a monomial of zero degree, i.e. $k_{j1}+ k_{j2}+\cdots +k_{jr}>0; k_{ji}\geq0$, for  all $j=1,...,N$. The special integral of a multidimensional Tarry's problem is defined as in [4, 8] as an integral

\[ \theta_{k}=\int_{-\infty}^{\infty}\int_{-\infty}^{\infty}\cdots \int_{-\infty}^{\infty}
\left|\int_{0}^{1}\cdots \int_{0}^{1}e^{2\pi i F(\bar{x})}d\bar{x}\right|^{2k}d\alpha_{1}d\alpha_{2}\cdots d\alpha_{N}.\]

\textit {Definition 1.} The  number $\gamma>0$  is called to be a convergence exponent of multi-dimensional Tarry's problem, if the special integral $\theta_{k}$  converges when $2k>\gamma$ , and diverges when $2k<\gamma$ .

For questions concerning the history of the problem one can refer to the works [7, 8, 10-15]. The authors of the works [4, 8] found an exact value of the convergence exponent of one-dimensional problem, and had observed an effect of variation of the value of a convergence exponent depending on a view (complete or incomplete) of the polynomial. A conjecture arose on variation of the value of the convergence exponent depending on the view of a concrete polynomial (see [4, p. 53). They proved estimation from above for the convergence exponent in multi-dimensional case for two types of polynomials. Both bounds based on the estimation of multiple trigonometrical integrals. Since the question is connected with the estimation of trigonometrical integrals, then the value of the convergence exponent seemed to be dependent on singularities connected with the given polynomial finally.

So, there was a firm confidence that in the multidimensional problem on convergence exponent for the special integral influence of singularities have a crucial importance (see [17, p.179, also [18]). However, there are errors in the works [17, 18], and some of lower bounds are not correct. Both methods of getting of lower bounds are the same and we indicate to the error in the work [17]. At the p. 183 (and further) of this work one made an exchange of variables in outer integral depending of parameters included into the inner integral. Despite that for fixed values of these parameters the Jacobian of this change is equal to 1, the correspondence of domains must be taken into attention (see [16, p. 255], or [19, p. 497]). Really, to consider these parameters as fixed we need in changing of order of integration. As a result of ignoring of dependence from parameters the lower estimations were got are incorrect. The same mistake was made in the the proof of the theorem 2 [18].

For the purposes of preliminary observation we consider the polynomial
\[F(x_{1} ,x_{2} ,...,x_{r} )=\sum _{j=1}^{m}\alpha _{j} (x_{1} x_{2} \cdots x_{r} )^{j}  ,\]
where $m$ is a natural number. The monomial map defined by this polynomial, i.e. the map
\[\varphi :(x_{1} ,...,x_{r} )\mapsto {}^{t} \left(\begin{array}{ccc} {(x_{1} x_{2} \cdots x_{r} )^{m} ,} & {\cdots ,} & {x_{1} x_{2} \cdots x_{r} } \end{array}\right)\]
has a Jacoby matrix:
\[\left(\begin{array}{ccc} {mx_{1}^{-1} (x_{1} x_{2} \cdots x_{r} )^{m} } & {\cdots } & {mx_{r}^{-1} (x_{1} x_{2} \cdots x_{r} )^{m} } \\ {(m-1)x_{1}^{-1} (x_{1} x_{2} \cdots x_{r} )^{m-1} } & {\cdots } & {(m-1)x_{r}^{-1} (x_{1} x_{2} \cdots x_{r} )^{m-1} } \\ {\vdots } & {\ddots } & {\vdots } \\ {x_{1}^{-1} (x_{1} x_{2} \cdots x_{r} )} & {\cdots } & {x_{r}^{-1} (x_{1} x_{2} \cdots x_{r} )} \end{array}\right).\]
It is easy to see that the rank of this matrix is equal to 1. So, each point of the unite cube $[0,1]^{r} $ is a singular point of the given map. However, if to consider the Cartesian product $[0,1]^{r} \times [0,1]^{r} \times \cdots \times [0,1]^{r} $ of  $k$ exemplars of the unite cube, and to consider the map (such maps arise in the estimations of the convergence exponent)
\begin{equation} \label{1}
\psi (\bar{x}_{1} ,\bar{x}_{2} ,...,\bar{x}_{k} )=\left(\varphi (\bar{x}_{1} )+\varphi (\bar{x}_{2} )+...+\varphi (\bar{x}_{k} )\right);\; \bar{x}_{j} =(x_{j1} ,...,x_{jr} )\in [0,1]^{r}
\end{equation}
then for $k\ge m$ it will be regular (i.e. its Jacoby matrix has a maximal rank) everywhere in the cube $[0,1]^{r} \times [0,1]^{r} \times \cdots \times [0,1]^{r} =[0,1]^{kr} $, except for points of some subset of a zero Jordan measure. Really, Jacoby matrix of the map (1) consists of $k$ blocks of a kind indicated above. If to take a submatrix of order $m$ of this matrix, containing more than two columns from any block, then the constructed minor will be equal to zero. At the same time, if we take the submatrix of order $m$ of the matrix composed by the first columns of each block, i.e., for example, of the matrix
\[\left(\begin{array}{ccc} {mx_{11}^{-1} (x_{11} x_{12} \cdots x_{1r} )^{m} } & {\cdots } & {mx_{k1}^{-1} (x_{k1} x_{k2} \cdots x_{kr} )^{m} } \\ {(m-1)x_{11}^{-1} (x_{11} x_{12} \cdots x_{1r} )^{m-1} } & {\cdots } & {(m-1)x_{k1}^{-1} (x_{k1} x_{k2} \cdots x_{kr} )^{m-1} } \\ {\vdots } & {\ddots } & {\vdots } \\ {x_{11}^{-1} (x_{11} x_{12} \cdots x_{1r} )^{} } & {\cdots } & {x_{k1}^{-1} (x_{k1} x_{k2} \cdots x_{kr} )} \end{array}\right)\]
then we will receive algebraic sum of $m!$ products which are parewisely different monomials.  Therefore, we get a polynomial of degree $[rm(m+1)-2m]/2$ of $mr$ variables with nonzero coefficients. Then (see [9]), this polynomial vanishes in the subset of the Cartesian product $[0,1]^{r} \times [0,1]^{r} \times \cdots \times [0,1]^{r} $ of a zero measure only, contained in by the union of finite number of hypersurfaces. Let's consider now the system connected with the given polynomial:
\[x_{11}^{m} ...x_{1r}^{m} +\cdots+x_{k1}^{m} ...x_{kr}^{m}-x_{k+1,1}^{m} ...x_{k+1,r}^{m} -\cdots -x_{2k,1}^{m} ...x_{2k,r}^{m} =0,\]
\[x_{11}^{m-1} ...x_{1r}^{m-1}+\cdots+x_{k1}^{m-1} ...x_{kr}^{m-1}-x_{k+1,1}^{m-1} ...x_{k+1,r}^{m-1} -\cdots -x_{2k,1}^{m-1} ...x_{2k,r}^{m-1} =0,\]
\[\cdots \quad \cdots \quad \cdots \]
\[x_{11}^{} ...x_{1r}^{} +\cdots+x_{k1}^{} ...x_{kr}^{}-x_{k+1,1}^{} ...x_{k+1,r}^{} -\cdots -x_{2k,1}^{} ...x_{2k,r}^{} =0.\]
The set of solutions of this system allows a shift. Namely, exchanging of the variables by the equality of a type $t_{j}=x_{j1} x_{j2} \cdots x_{jr} $ reduces the given polynomial to the one-dimensional complete polynomial case (see [4, 8]). If the set of solutions contains a vector $(t_{1} ,t_{2} ,...,t_{2k} )$, then it contains a new vector defined by the conditions
\[x'_{j1} x'_{j2} \cdots x'_{jr} =x_{j1} x_{j2} \cdots x_{jr} +a;\, j=1,...,2k,\]
for every real $a$. The equality $t=x_{1} x_{2} \cdots x_{r} $ does't define the variables $x_{1} ,x_{2} ,\cdots ,x_{r} $ one-valuedly. By this reason when $a$ varies continuously the all of $r$ variables $x_{1} ,x_{2} ,\cdots x_{r} $ varies unboundedly. For restore one-valudity we need to fix, for example, the all of the variables $x_{11} ,x_{12} ,\cdots ,x_{1r} $. Under these constraints we can estimate the surface integral by using of the results of the work [14]. Then,  as in one-dimensional case (see [8]), the integral converges when
\[2k>\frac{1}{r} \left(\frac{rm(m+1)}{2} +r\right)=\frac{m(m+1)}{2} +1.\]
Hence, singularities in the given problem does't affect an estimation found in [4] for the convergence exponent. Here it is important the fact that the set of solutions of the system considered above allows a shift. The similar situation arises in the general case.

\begin{center}
\textbf{2. Auxiliary results}
\end{center}

Let's consider a polynomial
\[F(x_{1} ,x_{2} ,...,x_{r} )=\sum _{j=1}^{N}\alpha _{j} \gamma _{j} (x_{1} ,x_{2} ,\cdots ,x_{r} ) ,\]
where $\gamma _{j} (\bar{x})=\gamma _{j} (x_{1} ,x_{2} ,\cdots ,x_{r} )$ are the monomials of a kind
\[\gamma _{j} (\bar{x})=x_{1}^{k_{j1} } x_{2}^{k_{j2} } \cdots x_{r}^{k_{jr} } .\]
Following by [4] we call the system
\[\gamma_{j}(\bar{x}_{1})+\cdots+\gamma_{j}(\bar{x}_{k})-\gamma_{j}(\bar{x}_{k+1})-\cdots-\gamma_{j}(\bar{x}_{2k})=0, j=1,...,N \]
to be the system connected with this polynomial for every natural $k$. The left sides of the equations of this system define a map  $\varphi ={}^{t} \left(\varphi _{1} ,\varphi _{2} ,...,\varphi _{N} \right)$, where
\[\varphi _{j} :(\bar{x}_{1} ,...,\bar{x}_{2k} )\mapsto
\left(\gamma_{j}(\bar{x}_{1})+\cdots+\gamma_{j}(\bar{x}_{k})-\gamma_{j}(\bar{x}_{k+1})-\cdots-\gamma_{j}(\bar{x}_{2k}\right)\]
 The Jacoby matrix looks (variables we assume non-zero) like:
\begin{equation} \label{2}
J_{0}=\left(\begin{array}{cccc} {k_{11} x_{11}^{-1} \gamma _{1} (\bar{x}_{1} )} & {k_{12} x_{12}^{-1} \gamma _{1} (\bar{x}_{1} )} & {\cdots } & {-k_{1r} x_{2k,r}^{-1} \gamma _{1} (\bar{x}_{2k} )} \\ {k_{21} x_{11}^{-1} \gamma _{2} (\bar{x}_{1} )} & {k_{22} x_{12}^{-1} \gamma _{2} (\bar{x}_{1} )} & {\cdots } & {-k_{2r} x_{2k,r}^{-1} \gamma _{2} (\bar{x}_{2k} )} \\ {\vdots } & {\vdots } & {\ddots } & {\vdots } \\ {k_{N1} x_{11}^{-1} \gamma _{N} (\bar{x}_{1} )} & {k_{N2} x_{12}^{-1} \gamma _{N} (\bar{x}_{1} )} & {\cdots } & {-k_{Nr} x_{2k,r}^{-1} \gamma _{N} (\bar{x}_{2k} )} \end{array}\right).
\end{equation}
The minor of order $r$ of this matrix, composed of the first $r$ columns, after of reducing by common factors of the elements of columns and rows of the determinant, leads to minors of a matrix
\begin{equation} \label{3}
\left(\begin{array}{cccc} {k_{11} } & {k_{12} } & {\cdots } & {k_{1r} } \\ {k_{21} } & {k_{22} } & {\cdots } & {k_{2r} } \\ {\vdots } & {\vdots } & {\ddots } & {\vdots } \\ {k_{N1} } & {k_{N2} } & {\cdots } & {k_{Nr} } \end{array}\right).
\end{equation}

Designate
\[S(\bar{x})=\left(\begin{array}{cccc} {k_{11} x_{1}^{-1} \gamma _{1} (\bar{x})} & {k_{12} x_{2}^{-1} \gamma _{1} (\bar{x})} & {\cdots } & {k_{1\rho} x_{\rho}^{-1} \gamma _{1} (\bar{x})} \\ {k_{21} x_{1}^{-1} \gamma _{2} (\bar{x})} & {k_{22} x_{2}^{-1} \gamma _{2} (\bar{x})} & {\cdots } & {k_{2\rho} x_{\rho}^{-1} \gamma _{2} (\bar{x})} \\ {\vdots } & {\vdots } & {\ddots } & {\vdots } \\ {k_{N1} x_{1}^{-1} \gamma _{N} (\bar{x})} & {k_{N2} x_{2}^{-1} \gamma _{1} (\bar{x})} & {\cdots } & {k_{N\rho} x_{\rho}^{-1} \gamma _{1} (\bar{x})} \end{array}\right).\]
and
\[K(\bar{x})=\left(\begin{array}{cccc} {k_{11} \gamma _{1} (\bar{x})} & {k_{12} \gamma _{1} (\bar{x})} & {\cdots } & {k_{1\rho } \gamma _{1} (\bar{x})} \\ {k_{21} \gamma _{2} (\bar{x})} & {k_{22} \gamma _{2} (\bar{x})} & {\cdots } & {k_{2\rho } \gamma _{2} (\bar{x})} \\ {\vdots } & {\vdots } & {\ddots } & {\vdots } \\ {k_{N1} \gamma _{N} (\bar{x})} & {k_{N2} \gamma _{N} (\bar{x})} & {\cdots } & {k_{N\rho } \gamma _{N} (\bar{x})} \end{array}\right).\]
for each  $\bar{x}=(x_{1},x_{2},...,x_{r})$. Then the Jacoby matrix (2) can be written as a block-matrix
 \[\varphi'=(S(\bar{x}_{1}),...,S(\bar{x}_{k}),  -S(\bar{x}_{k+1}),...,-S(\bar{x}_{2k})).\]
Arranging the entries of every column consequently in a row, take the transposed Jacoby matrix of the obtained system of monomials, designating the received matrix as  $\varphi''$. Continuing this procedure we get the sequence of matrices  $\varphi^{(j)}, j=1,..,m$, where $m=deg F$ . Let $W_{j}$  be the set of points of $[0,1]^{2k}$ at which $\Phi_{j}=det(\varphi^{(j)}\cdot{}^{t}\varphi^{(j)})=0, j=0,1,...,m   $ . The set of singular points is a union $W=\bigcup_{j=1}^{m}W_{j}$ . Following two lemmas show that the algebraic set $W$ has a zero Jordan measure (for the proof of our theorems below we need in a more detailed structure of algebraic sets (see [9])).

To formulate our lemmas we need in a definition.

\textit{Definition 1.} The matrix $A$ of a rank $\rho$  we will call to be the matrix of the structure $(\rho_{1},...,\rho_{q} )$ if it can be represented as a block-matrix
\[A=\left(\begin{array}{c}{A_{1}}\\{\cdots}\\{A_{q}} \end{array}\right),\]
and: 1) $rank A_{j}=\rho_{j}$; 2) if to the matrix $A_{j}$ to join the first row of the matrix $A_{j+1}$ then its rank remains unchanged, but if to delete its last row then the rank of the matrix $A_{j}$ stands less.

Let $\mathbf{N_{0}}$ be a set of non-zero integral numbers. In this set operation of addition is defined and an algebra $\langle \mathbf{N_{0}}, + \rangle$ is a commutative monoid.	Let us consider an algebra $\langle \mathbf{N_{0}}^{r}, + \rangle$.

\textit{Definition 2.} A strict order $\prec$ we call to be a normally-lexicographic (n.-l.) order, if following relations are satisfied:
\[(k_{j1},...,k_{jr})\prec (k_{i1},...,k_{ir}) \Leftrightarrow k_{j1}+\cdots +k_{jr}<k_{i1}+\cdots +k_{ir}\bigvee\]
\[\bigvee (k_{j1}+\cdots +k_{jr}=k_{i1}+\cdots +k_{ir}\bigwedge (k_{j1},...,k_{jr})\prec_{1} (k_{i1},...,k_{ir})),\]
where $\prec_{1}$ denotes the usual lexicographic order.

\textit {Example.} A set of pares
\[M=\{(0,1),(1,0),(0,2),(1,1),(2,0),(0,3),(1,2),(2,1),(3,0)\} \]
is an increasing n.-l. ordered set.

\textit{Proposition 1.} The algebraic system $\langle \mathbf{N_{0}^{r}}, +, \prec \rangle$ is an ordered monoid in which the strict and linear order $\prec$ satisfies the condition:
\[ (\forall a,b,c,d\in \mathbf{N_{0}^{r}})(a\prec c\bigwedge b\prec d\Rightarrow a+c\prec b+d).\]

\textit{Proof.} Let $a=(a_{1},...,a_{r}), b=(b_{1},...,b_{r}), c=(c_{1},...,c_{r}), d=(d_{1},...,d_{r})$. Consider two cases:

1) $a_{1}+\cdots +a_{r}<c_{1}+\cdots +c_{r} \bigvee b_{1}+\cdots +b_{r}<d_{1}+\cdots +d_{r};$

2) $a_{1}+\cdots +a_{r}=c_{1}+\cdots +c_{r} \bigwedge b_{1}+\cdots +b_{r}=d_{1}+\cdots +d_{r}.$
In the first case we have:
\[(a_{1}+\cdots +a_{r})+(b_{1}+\cdots +b_{r})<(c_{1}+\cdots +c_{r})+(d_{1}+\cdots +d_{r}),\]
because from the relation $b\prec d$, in consent with the definition of the order, we deduce that if $a_{1}+\cdots +a_{r}<c_{1}+\cdots +c_{r}$ then
\[b_{1}+\cdots +b_{r}=d_{1}+\cdots +d_{r}\bigvee b_{1}+\cdots +b_{r}<d_{1}+\cdots +d_{r}.\]
In both cases we get $a+c\prec b+d$.

Consider now the case 2). Let the conditions $a_{1}+\cdots +a_{r}=c_{1}+\cdots +c_{r}\bigvee b_{1}+\cdots +b_{r}=d_{1}+\cdots +d_{r}.$ Then
\[(a_{1}+\cdots +a_{r})+(b_{1}+\cdots +b_{r})=(c_{1}+\cdots +c_{r})+(d_{1}+\cdots +d_{r}),\]
and the required result follows from the properties of lexicographic order. The proposition is proven.

\textit{Definition 3.} Let $F(x_{1},...,x_{r})\in\mathbf{R}[x_{1},...,x_{r}]$ be any polynomial.
We say that the monomials of the polynomial are n.-l. ordered if the set $M$ of exponents' corteges of monomials is a n.-l. ordered set. The monomial corresponding to the high member of the set $M$ we call to be n.-l. high member of the polynomial.

The following proposition is an immediate consequence of the proposition 1.

\textit{Proposition 2.} The n.-l. high member of the product of two polynomials is equal to the product of high members of given polynomials.

\textit{Proposition 3.} Let rows of the matrix (3) be n.-l. ordered, have a rank $\rho$ and the structure $(\rho_{1},...,\rho_{q})$. Then the matrix
\[(K(\bar{x}_{1}), K(\bar{x}_{2}),...,K(\bar{x}_{q}))\]
has a rank $N$ for all $(\bar{x}_{1}, \bar{x}_{2},...,\bar{x}_{q})$ with exception of their set of zero Jordan measure.

\noindent \textit{ Proof. } It is clear that $\rho_{1} \leq \rho,...,\rho_{q} \leq \rho$ and $\rho_{1} +\cdots +\rho_{q}=N$ (this follows from that that the matrix (3) has not zero row). Note that if $\rho=N$ then the statement of the proposal is true for the matrix $K(\bar{x}_{1})$ already. So, we suppose that $\rho<N$. Since the matrix has a structure $(\rho_{1},...,\rho_{q})$, then the matrix $K(\bar{x}_{2})^{\rho_{2}}$ consisted in the rows of the matrix $K(\bar{x}_{2})$ having the rank $\rho_{2}$. Take a base minor and denote it as $K_{2}$. Consider now the block-matrix
\[Y=\left( \begin{array}{cc}{K_{1}}&{K'_{2}}\\{K'_{1}}&{K_{2}} \end{array}\right),\]
where the blocks $K'_{1}$ and $K'_{1}$ consisted of the entries of the matrix $(K(\bar{x}_{1}) K(\bar{x}_{2}))$ standing on the intersection of the columns and rows, where the base minors $K_{1}$ and $K_{2}$ are placed. Take the expanding of the determinant of the matrix $Y$ in cofactors of  minors of the columns on which the base minor $K_{1}$ is placed. Into the expanding is entered the product $K_{1}K_{2}$. The minor $K_{2}$ is a monomial of variables of the vector $\bar{x}_{2}$ with non-zero coefficient. Moreover, in consent with the proposition 2, this monomial is a high monomial which contained by the matrix
\[\left(\begin{array}{c}{K'_{2}}\\{K_{2}} \end{array}\right),\]
and by this reason will be a high member of the polynomial got after of expanding of the determinant of the matrix $Y$. Consequently, this determinant is distinct from zero everywhere with exception of the points a set of zero Jordan measure. It is obvious that we can repeat taking a new block matrix with the matrix on left top corner. We arrive at the proof of the proposal 3, continuing by this way.

Consider separately the case of two-dimensional problem, and take the following polynomial ($n\ge m$) with n. -l. ordered members:
\[ F(x,y)=\sum_{i=1}^{n}\sum_{j=1}^{m} \alpha_{ij}x^{i}y^{j} .\]
Denote by $M_{n}$ the transposed matrix of exponents i. e.

\[M_{n}=\left(\begin{array}{ccccccc} {0}&{1}&{0}&{1}&{2}&{\cdots}&{n}\\{1}&{0}&{2}&{1}&{0}&{\cdots}&{m}\end{array}\right). \]

\noindent  \textbf{\textit{ Lemma 1. }} If  the number of columns of this matrix is an even number then it has a structure $(2,2,\dots,2)$, if else, it has a structure $(2,2,\dots,2,1)$.

	\noindent \textit{ Proof.} Consider two cases: 1) the number of columns of the matrix $M_{n}$ is an even number; 2) the number of columns of the matrix $M_{n}$ is an odd number.

1) Dissect the matrix $M_{n}$  into blocks of an order 2. We will call a block to be the block of a stage $k$ if the sum of entries of each its columns is equal $k$. Every of remaining blocks we will call to be mixed. For example the blocks $\left(\begin{array}{cc}{0}&{1}\\{4}&{3}\end{array}\right)$, $\left(\begin{array}{cc} {2}&{3}\\{2}&{1} \end{array} \right)$ belong to the stage 4, but the block $\left(\begin{array}{cc} {4}&{0}\\{0}&{5} \end{array} \right)$ is a mixed block. All of mixed blocks have a view $\left(\begin{array}{cc} {k}&{0}\\{0}&{k+1}\end{array} \right) $ (if $k+1\le m$), or a view $\left(\begin{array}{cc} {k}&{k-m+1}\\{0}&{m}\end{array} \right) $ (if $n+1\ge k+1 >m$), or $\left(\begin{array}{cc} {n}&{k+1-m}\\{k-n}&{m}\end{array} \right) $ (if $m+n-1\ge k>n$). Since
\[ \det \left(\begin{array}{cc} {n}&{n-m+1+a}\\{a}&{m}\end{array}\right) =nm-an+am-a-a^{2}\ge n(m-a)>0,\]
($m\ge 1+a$) then every mixed block is non-singular. Since a block of a stage $k$ has a view  $\left(\begin{array}{cc} {a}&{a+1}\\{b}&{b-1}\end{array} \right) $, with $a+b=k, a\ge 0$, then they are non-singular also. So, the lemma is proven for the case 1).
	
In the case 2) all of said above remains without changing for the first $n-1$ columns of the matrix $M_{n}$. So the lemma 1 true in this case also. Proof of the lemma 1 is finished.

\textit{Remark.} In the work [13] the matrix of exponents (3) satisfies the conditions of the lemma 1. So, the statement of the lemma 1 is true, which was accepted without proof in [13].
 	
In the works [13-14] the matrices $A_{j}; j=1, 2,...$ were introduced by taking of transposed Jacoby matrices of monomials. So, every matrix  $A_{j}; j\geq 1$  is composed by blocks, and every block is a transposed Jacoby matrix of the system of all partial derivatives of some monomial.

\noindent  \textbf{\textit{ Lemma 2. }}	Let a monomial $ax_{1}^{k_{1}} x_{2}^{k_{2}}\cdots x_{r}^{k_{r}}$  have non-zero multiplier $a$ and positive exponents $k_{1}, k_{2},..., k_{r}>0$, where $r>1$. Then Gram determinant of the system of coordinate-functions of the gradient of this monomial has a maximal rank for all $\bar{x}$  with exception of the points of some subset of zero Jordan measure.
	
\noindent \textit {Proof.} The gradient of a given monomial (we put $a=1$ ) is an column with coordinates

\[\left(\begin{array}{c} {k_{1} x_{1}^{k_{1}-1} x_{2}^{k_{2}}\cdots x_{r}^{k_{r}}} \\ {k_{2} x_{1}^{k_{1}} x_{2}^{k_{2}-1}\cdots x_{r}^{k_{r}}}\\ {\vdots } \\ {k_{r} x_{1}^{k_{1}} x_{2}^{k_{2}}\cdots x_{r}^{k_{r}-1}} \end{array}\right).\]
The Gram determinant of this system of functions has a view $A\cdot ^{t}A$ , and the matrix $A$  is defined as follows:

\[A=\left(\begin{array}{ccc} {k_{1}(k_{1}-1)x_{1}^{k_{1}-2} x_{2}^{k_{2}}\cdots x_{r}^{k_{r}}} & {\cdots } & {k_{r}k_{1} x_{1}^{k_{1}-1} x_{2}^{k_{2}}\cdots x_{r}^{k_{r}-1}} \\ {k_{1} k_{2} x_{1}^{k_{1}-1} x_{2}^{k_{2}-1}\cdots x_{r}^{k_{r}}} & {\cdots } & {k_{r}k_{2} x_{1}^{k_{1}} x_{2}^{k_{2}-1}\cdots x_{r}^{k_{r}-1}} \\ {\vdots } & {\ddots } & {\vdots } \\ {k_{1} k_{r} x_{1}^{k_{1}-1} x_{2}^{k_{2}}\cdots x_{r}^{k_{r}-1}} &{\cdots} & {k_{r}(k_{r}-1) x_{1}^{k_{1}} x_{2}^{k_{2}}\cdots x_{r}^{k_{r}-2}} \end{array}\right).\]
(note that if $k_{j}=1$  then corresponding monomial induces a zero entry after of repeated differentiation, as it is seen from the view of the matrix $A$). Taking the determinant of this matrix, it easy to establish that we can reduce the rows of this matrix correspondingly by the factors $k_{1}, k_{2},...,k_{r}$ . Then we may multiply the rows and columns by  $x_{1}, x_{2},...,x_{r}$, and we get the determinant

\[ \left | \begin{array}{cccc} {k_{1} -1} & {k_{2}}& {\cdots } & {k_{r}} \\ {k_{1}} & {k_{2}-1} & {\cdots } & {k_{r}} \\ {\vdots } & {\vdots } & {\ddots } & {\vdots } \\ {k_{1}} & {k_{2}} & {\cdots } & {k_{r}-1} \end{array}\right|,\]
after of reducing the columns by the monomial $x_{1}^{k_{1}} x_{2}^{k_{2}}\cdots x_{r}^{k_{r}}$  . Since the got determinant equals
\[ (-1)^{r}(1-k_{1}-\cdots -k_{r})\neq 0,\]
then the proof of the lemma 2 is ended.

\textit{Consequence.} Let the polynomial $F(\bar{x})$  be defined as above, and its senior form contains all of independent variables. Then the all of the matrices $\varphi^{(j)}=A_{j}(\bar{x}); 1 \leq j \leq m-1$  have a maximal rank for all points of the cube $[0,1]^{kr}$, with exception for the points of a subset of zero Jourdan measure.

\textit{Proof.}  When $j=1$  the statement of the consequence follows from the lemma 1. If $j=2$  then the matrix $A_{1}$  consisted of blocks every of which being transposed Jacoby matrix of the system of partial derivatives of every monomial of the polynomial $\gamma_{j}(\bar{x}_{1})+\cdots+\gamma_{j}(\bar{x}_{k})$. If all of variables enter into one monomial then the statement of the consequence is true due to theorem 2. Let now the variables $x_{11},...,x_{1r}$  contained by several monomials, for example by two monomials $x_{11}^{k_{1}}\cdots x_{1p}^{k_{p}}$  and  $x_{1l}^{q_{l}}\cdots x_{1r}^{q_{r}}$ of the considered above polynomial. The case $p=r-1$ is trivial, so we assume that $p<r-1$. Then by the lemma 2 the matrix $A_{2}$  contains two non-singular submatrices of order $p$  and   $r-l+1 (p+r-l+1\ge r )$. The rows on which these blocks are placed contain all of rows containing partial derivatives with respect to the variables $x_{11},...,x_{1r}$ . By the lemma 1 the rows corresponding to the variables $x_{1m+1},...,x_{1r}$  being linearly independent may contain the rows containing some of rows of the second block. But the rows containing the partial derivatives of the monomial $x_{1l}^{q_{l}}\cdots x_{1r}^{q_{r}}$  with respect to the variables  $x_{1p+1},...,x_{1r}$ doesn't contain partial derivatives of the monomial $x_{11}^{k_{1}}\cdots x_{1p}^{k_{p}}$  with respect to the variables $x_{1p+1},...,x_{1r}$ . So, the second block contain a submatrix being a transposed Jacoby matrix of a gradient of the monomial $x_{1p+1}^{q_{p+1}}\cdots x_{1r}^{q_{r}}$  . In concent with the lemma 1 this matrix is non-singular. Construct the minor consisting of two non-singular blocks placed on the diagonal: the first block and the constructed block of order $r-p$ . So, the matrix $A_{2}$  has a maximal rank.
By an analogy we can prove the statement of the consequence for every matrix $A_{j}$ . The consequence is proven.

\textbf{\textit{ Lemma 3. }}Let the polynomial $F(\bar{x})$ of degree $m$ be defined by the equality (2), and the matrix (3) have a rank $\rho \ge 1$. Then, for $k\rho \ge N$ and $k\ge q$ (when $q$ is defined by the conditions of the proposition 3) the following formula holds
\[\theta _{k} =(2\pi )^{N} \int _{\Pi }\frac{ds}{\sqrt{G_{0} } }  ,\]
and the integral on the right part of this equality is defined as a limit:
\[\mathop{\lim }\limits_{\eta \to 0} (2\pi )^{N} \int _{\Pi \bigcap D_{\lambda } \times D_{\lambda } }\frac{ds}{\sqrt{G_{0} } }  ;\]
here $D_{\lambda } =\left\{\left(\bar{x}_{1} ,...,\bar{x}_{k} \right)|min_{1\leq j \leq m-1}\Phi_{j} \ge \lambda \right\}$, $G_{0}=det(J_{0}\cdot ^{t}J_{0})$ ($J_{0}$ defined by (2)), $\Pi $ is an algebraic set defined by the system of the equations
\begin{equation} \label{8}
\gamma _{j} (\bar{x}_{1} )+\cdots +\gamma _{j} (\bar{x}_{k} )-\gamma _{j} (\bar{x}_{k+1} )-\cdots -\gamma _{j} (\bar{x}_{2k} )=0,\, \, \, j=1,...,N
\end{equation}
in $2kr$-dimensional unite cube.

\textit{Note.} This lemma is a correction of our result from [14] where the case $\rho=r$  was considered only without specification. At the same time there was used theorem of Sard the application of which demands some conditions. Here we bypass the use of the theorem of Sard using our lemmas 1 and 2 to isolate singularities by considering improper integrals.

\textit{Proof.} If $\rho =r$ then proof of the lemma 3 can be conduct as in [15]. So, we can assume, that $\rho <r$. It is enough to consider the case when $\rho =r-1$, because other cases can be considered similarly. Let one of columns of the matrix, say the last column, be a linear combination of the first $r-1$ columns, i. e. there are constants $\beta _{1} ,...,\beta _{r-1} \in \mathbf{R}$ such that
\[\beta _{1} k_{j1} +\cdots +\beta _{r-1} k_{jr-1} =k_{jr} ,\, \, j=1,...,N.\]
So,
\[\gamma _{j} (\bar{x})=\left(x_{1} x_{r}^{\beta _{1} } \right)^{k_{j1} } \left(x_{2} x_{r}^{\beta _{2} } \right)^{k_{j2} } \cdots \left(x_{r-1} x_{r}^{\beta _{r-1} } \right)^{k_{jr-1} } =\eta _{j} (\bar{y}),\]
and, consequently, after of changing of variables of a kind $y_{{\rm i}} =x_{i} x_{r}^{\beta _{i} } $, we receive the polynomial
\[G(y_{1} ,y_{2} ,...,y_{r-1} )=\sum _{j=1}^{N}\alpha _{j} \delta _{j} (y_{1} ,y_{2} ,...,y_{r-1} ) \]
of  $r-1$ variables, where
\[\delta _{j} (y_{1} ,y_{2} ,...,y_{r-1} )= y_{1}^{k_{j1} } \cdots y_{r{\rm -1}}^{k_{jr-1} } .\]
For the polynomial $G(y_{1} ,y_{2} ,...,y_{r-1} )$ the matrix (3) has the maximal rank. Hence, the general case, as well as above, after of some transformations is reduced to the basic case when the matrix (3) has the maximal rank. We have:
\[\int _{0}^{1}\cdots \int _{0}^{1}\int _{0}^{1}e^{2\pi iF(x_{1} ,...,x_{r} )}  dx_{1} \cdots dx_{r}=\]
\[=\int _{0}^{1}x_{r}^{-\beta } dx_{r} \int _{0}^{x_{r}^{\beta _{1} } }\cdots \int _{0}^{x_{r}^{\beta _{r} } }e^{2\pi iG(y_{1} ,...,y_{:-1} )}  dy_{1} \cdots dy_{r-1}   =\]
\begin{equation} \label{9}
=\int _{0}^{\infty }\cdots \int _{0}^{\infty }\left(\int _{\chi _{0} (y_{1} ,...,y_{r-1} )}^{\chi _{1} (y_{1} ,...,y_{r-1} )}x_{r}^{-\beta } dx_{r}  \right)e^{2\pi iG(y_{1} ,...,y_{:-1} )} dy_{1} \cdots dy_{r-1}   ,
\end{equation}
where we have designated $\beta =\beta _{1} +\beta _{2} +\cdots +\beta _{r-1} $, and the functions $\chi _{0} (y_{1} ,...,y_{r-1} )$ and $\chi _{1} (y_{1} ,...,y_{r-1} )$ are defined by the condition that the all of inequalities $0<x_{i} \le 1$ are carried out simultaneously. Designating
\[\varphi (y_{1} ,...,y_{r-1} )=\int _{\chi _{0} (y_{1} ,...,y_{r-1} )}^{\chi _{1} (y_{1} ,...,y_{r-1} )}x_{r}^{-\beta } dx_{r}  ,\]
we can rewrite (9) as a follows
\[\int _{0}^{1}\cdots \int _{0}^{1}\int _{0}^{1}e^{2\pi iF(x_{1} ,...,x_{r} )}  dx_{1} \cdots dx_{r} =\]
\[=\int _{0}^{\infty }\cdots \int _{0}^{\infty }\varphi (y_{1} ,...,y_{r-1} )e^{2\pi iG(y_{1} ,...,y_{r-1} )} dy_{1} \cdots dy_{r-1}   .\]
As well as in the basic case, considering the $k$- th power of the trigonometrical integral and spending the reasonings applied above, we will represent it as a Fourier transformation:
\[\left(\int _{0}^{1}\cdots \int _{0}^{1}\int _{0}^{1}e^{2\pi iF(x_{1} ,...,x_{r} )}  dx_{1} \cdots dx_{r}   \right)^{k} =\]
\[=\int _{0}^{\infty }\cdots \int _{0}^{\infty }\cdots \cdots \int _{0}^{\infty }\cdots \int _{0}^{\infty }\varphi (y_{11} ,...,y_{1\, \, r-1} )\cdots \varphi (y_{k1} ,...,y_{k\, \, r-1} )  \times   \]
\[\times e^{2\pi iG(y_{11} ,...,y_{1r-1} )} \cdots e^{2\pi iG(y_{k1} ,...,y_{k\, \, r-1} )} dy_{11} \cdots dy_{k\, \, r-1} =\]
\[=\int _{}^{}\cdots \int _{}^{}\cdots \cdots \int _{}^{}\cdots \int _{}^{}\varphi (y_{11} ,...,y_{1\, \, r-1} )\cdots \varphi (y_{k1} ,...,y_{k\, \, r-1} )\times     \]
\[\times e^{2\pi i\left\{\alpha _{1} \sum _{i=1}^{k}\eta _{{\rm 1}} (\bar{y}_{i} )+\cdots +\alpha _{N} \sum _{i=1}^{k}\eta _{{\rm N}} (\bar{y}_{i} )  \right\}} dy_{11} \cdots dy_{k\, \, r-1} =\]
\[=\int _{0}^{k}\cdots \int _{0}^{k}e^{2\pi i(\alpha _{1} u_{1} +\cdots +\alpha _{N} u_{N} )}   \int _{\Pi (\bar{u})}\varphi (y_{11} ,...,y_{1\, \, r-1} )\cdots \varphi (y_{k1} ,...,y_{k\, \, r-1} )\frac{ds}{\sqrt{G} }  ;\]
here $\bar{u}=(u_{1} ,...,u_{N} )$, and $\Pi (\bar{u})=\Pi (u_{1} ,...,u_{N} )$ designates a surface defined in $[0,\infty ]^{N} $ by the system of  equations
\[\sum _{i=1}^{k}\eta _{j} (\bar{y}_{i} ) =u_{j} ,\, \, j=1,...,N, \bar{u}=(u_{1} ,...,u_{N} )\in [0,k]^{N} ,\]
and $G$ denotes the Gram determinant of gradients of the functions standing on the left parts of the system's equations. For fixed $\bar{u}$ the surface $\Pi (\bar{u})$ is an image of one-valued mapping $y_{im} =x_{im} x_{rm}^{\beta _{i} } ,\, \, y_{ir} =x_{ir} ,\, m=1,...,r-1;\, i=1,...,k$ of the surface $\Pi '(\bar{u})$ in the unite cube $[0,1]^{kr} $.

Let's consider the transformation $y_{m} =x_{m} x_{r}^{\beta _{m} } ,\, \, m=1,...,r-1;\, y_{r} =x_{r} $ in the cube $[0,1]^{r}$. It has a Jacoby matrix
\[J=\left(\begin{array}{cccc} {x_{r}^{\beta _{1} } } & {{\rm 0}} & {\cdots } & {\beta _{1} x_{1} x_{r}^{\beta _{1} -1} } \\ {0} & {x_{r}^{\beta _{2} } } & {\cdots } & {\beta _{2} x_{2} x_{r}^{\beta _{2} -1} } \\ {\vdots } & {\vdots } & {\ddots } & {\vdots } \\ {0} & {0} & {\cdots } & {1} \end{array}\right),\]
with the determinant $x_{r}^{\beta -\beta _{r} } $. Hence, this mapping is one valued for $x_{r}^{} \ne 0$, which we suppose to be true. Then, the system of functions $\eta {}_{j} (\bar{y}),\, j=1,...,N$ will have the Jacoby matrix
\[\left(\begin{array}{cccc} {k_{11} y_{1}^{k_{11} -1} \cdots y_{r-1}^{k_{11} } } & {k_{12} y_{1}^{k_{11} } \cdots y_{r-1}^{k_{11} } } & {\cdots } & {k_{1r-1} y_{1}^{k_{11} } \cdots y_{r-1}^{k_{1r-1} -1} } \\ {k_{21} y_{1}^{k_{21} -1} \cdots y_{r-1}^{k_{21} } } & {k_{22} y_{1}^{k_{21} } \cdots y_{r-1}^{k_{21} } } & {\cdots } & {k_{2r-1} y_{1}^{k_{21} } \cdots y_{r-1}^{k_{2r-1} -1} } \\ {\vdots } & {\vdots } & {\ddots } & {\vdots } \\ {k_{N1} y_{1}^{k_{N1} -1} \cdots y_{r-1}^{k_{21} } } & {k_{N2} y_{1}^{k_{N1} } \cdots y_{r-1}^{k_{N1} } } & {\cdots } & {k_{Nr-1} y_{1}^{k_{N1} } \cdots y_{r-1}^{k_{Nr-1} -1} } \end{array}\right)\times \]
\[\times \left(\begin{array}{ccccc} {x_{r}^{\beta _{1} } } & {0} & {\cdots } & {0} & {\beta _{1} x_{1} x_{r}^{\beta _{1} -1} } \\ {0} & {x_{r}^{\beta _{2} } } & {\cdots } & {0} & {\beta _{2} x_{2} x_{r}^{\beta _{2} -1} } \\ {\vdots } & {\vdots } & {\ddots } & {\vdots } & {\vdots } \\ {0} & {0} & {\cdots } & {x_{r}^{\beta _{r-1} } } & {\beta _{r-1} x_{r-1} x_{r}^{\beta _{r-1} -1} } \end{array}\right).\]
It is easy to see that the entry standing on $i$-th line and $j$- th column is equal:
\[k_{ij} y_{1}^{k_{i1} } \cdots y_{j}^{k_{ij} -1} \cdots y_{r-1}^{k_{ir-1} } x_{r}^{\beta _{j} } =\]
\[=k_{ij} x_{1}^{k_{i1} } \cdots x_{j}^{k_{ij} -1} \cdots x_{r-1}^{k_{i1} } x_{r}^{k_{i1} \beta _{1} +\cdots +k_{ij} \beta _{j} -\beta _{j} +\cdots +k_{ir-1} \beta _{r-1} +\beta _{j} } =\]
\[=k_{ij} x_{1}^{k_{i1} } \cdots x_{j}^{k_{ij} -1} \cdots x_{r-1}^{k_{i1} } x_{r}^{k_{ir} } ,\]
if $j<r$. Consider now the case $j=r$. Multiplying the $i$-th line to the $j$-th column, we receive:
\[\sum _{j=1}^{r-1}k_{ij}  y_{1}^{k_{i1} } \cdots y_{j}^{k_{ij} -1} \cdots y_{r-1}^{k_{ir-1} } \beta _{j} x_{j} x_{r}^{\beta _{j} -1} =\]
\[=\sum _{j=1}^{r-1}k_{ij} \beta _{j}  x_{1}^{k_{i1} } \cdots x_{j}^{k_{ij} } \cdots x_{r-1}^{k_{i1} } x_{r}^{k_{i1} \beta _{1} +\cdots +k_{ij} \beta _{j} +\cdots +k_{ir-1} \beta _{r-1} -1} =\]
\[=k_{ir} x_{1}^{k_{i1} } \cdots x_{j}^{k_{ij} } \cdots x_{r}^{k_{ir} -1} .\]
Obviously, the Gram determinant $G$ will be transformed into $G_{0} $. Applying the lemma 2, of [12], we receive:
\[\int _{\Pi (\bar{u})}\varphi (y_{11} ,...,y_{1\, \, r-1} )\cdots \varphi (y_{k1} ,...,y_{k\, \, r-1} )\frac{ds}{\sqrt{G} }  =\]
\[=\int _{\Pi '(\bar{u})}\varphi (y_{11} ,...,y_{1\, \, r-1} )\cdots \varphi (y_{k1} ,...,y_{k\, \, r-1} )(x_{1r} \cdots x_{2k\, r} )^{\beta } \frac{d\sigma }{\sqrt{G_{0} } }  =\int _{\Pi '(\bar{u})}\frac{d\sigma }{\sqrt{G_{0} } }  .\]
So, we in accuracy have received the equality (see [14, p. 78]):
\[\theta _{k} =(2\pi )^{N} \int _{0}^{k}\cdots \int _{0}^{k}\left(\int _{\Pi '(\bar{u})}\frac{ds}{\sqrt{G_{0} } }  \right)  ^{2} du_{1} \cdots du_{N} .\]

The deduction of the formula of the lemma 3 from this relation can now be spent by a known scheme (see [13, 15]).
The lemma 3 is proven.

\begin{center}
\noindent \textbf{3. Basic results}
\end{center}

\textbf{Theorem 1.} If $k$ is a natural number such that $k\rho <N$ then $\theta _{k} $ diverges.

\noindent  \textit{ Proof. }We have:
\[\left(\int _{0}^{1}\cdots \int _{0}^{1}\int _{0}^{1}e^{2\pi iF(x_{1} ,...,x_{r} )}  dx_{1} \cdots dx_{r}   \right)^{k} =\]
\[=\int _{0}^{k}\cdots \int _{0}^{k}e^{2\pi i(\alpha _{1} u_{1} +\cdots +\alpha _{N} u_{N} )}   \int _{\Pi (\bar{u})}\varphi (y_{11} ,...,y_{1\, \, r-1} )\cdots \varphi (y_{k1} ,...,y_{k\, \, r-1} )\frac{ds}{\sqrt{G} }  ;\]
here $\bar{u}=(u_{1} ,...,u_{N} )$, $\Pi (\bar{u})=\Pi (u_{1} ,...,u_{N} )$ designates a surface defined in unite cube $[0,1]^{N} $ by the system of the equations
\[\sum _{i=1}^{k}\eta _{j} (\bar{y}_{i} ) =u_{j} ,\, \, j=1,...,N, \bar{u}=(u_{1} ,...,u_{N} )\in [0,k]^{N} ,\]
and $G$ is a Gram determinant of gradients of the functions standing on the left parts of the equations of the system. According to conditions of the theorem, the number of variables is less than number of equations. We will consider the first $k\rho $ equations of the system. The Jacoby matrix of the functions standing on the left parts of the considered system, owing to the lemma 2, has the maximal rank everywhere with exception of the points of a subset of zero Jourdan measure. Then the system of the considered equations defines, in a suitable neighbourhood of each point on the surface which doesn't belong to the exceptional set, every $\bar{y}_{j} $ as a vector-function of variables $u_{1} ,...,u_{k\rho } $. For $j>k\rho $ $u_{j} $ is not independent, and is a function of independent variables $u_{1} ,...,u_{k\rho } $ defined from other equations after of substitution of values of found variables. By the Parseval's formula
\[\int _{-\infty }^{\infty }\cdots \int _{-\infty }^{\infty }\left|\int _{0}^{1}\cdots \int _{0}^{1}e^{2\pi iF(\bar{x})} d\bar{x}  \right|  ^{2k} d\alpha _{1} ...d\alpha _{k\rho } =\]
\[=\int _{0}^{k}\cdots \int _{0}^{k}\left(\int _{\Pi (\bar{u})}\varphi (y_{11} ,...,y_{1\, \rho } )\cdots \varphi (y_{k1} ,...,y_{k\, \rho } )\frac{ds}{\sqrt{G} }  \right)^{2} du_{1} \cdots du_{k\rho }   .\]
All the multiple integrals are convergent or not simultaneously. The expression under the sign of the integral, by the told above, does't vanishes. Therefore, the last integral is non-zero and the further integration over $\alpha _{kr+1} ...\alpha _{N} $ through the interval $(-\infty ,+\infty )$ gives a divergent integral:
\[\int _{-\infty }^{\infty }\cdots \int _{-\infty }^{\infty }d\alpha _{k\rho +1} \cdots d\alpha _{N}\times\]
\[\times \int _{0}^{k}\cdots \int _{0}^{k}\left(\int _{\Pi (\bar{u})}\varphi (y_{11} ,...,y_{1\, \rho } )\cdots \varphi (y_{k1} ,...,y_{k\, \rho } )\frac{ds}{\sqrt{G} }  \right)^{2} du_{1} \cdots du_{k\rho }  .\]
The theorem 1 is proved.

 \textit{Definition 4. }We call the polynomial $F(\bar{x})$ to be $v$-complete if after of some permutation of variables $x_{1} ,x_{2} ,...,x_{r} $ this polynomial contains with any monomial $x_{1}^{k_{1} } ...x_{v}^{k_{v} } x_{v+1}^{k_{v+1} } ...x_{r}^{k_{r} } $ the all of monomials of a view $x_{1}^{k'_{1} } ...x_{v}^{k'_{v} } x_{v+1}^{k_{v+1} } ...x_{r}^{k_{r} } $ with
\[ 0\leq k'_{1} \leq k_{1},...,0\leq k'_{v} \leq k_{v}, 0\leq k'_{1} + \cdots+k'_{v} \le k_{v},\]
\[1\leq k'_{1} + \cdots+k'_{v}+k_{v+1}+\cdots+k_{r} \]

\textbf{Theorem 2. }Let the polynomial (2) be $v$-complete ($0\le v\le r$), and the matrix of exponents (3) has a rank $\rho $, $1\le \rho \le r$, is a matrix of the structure $(\rho_{1},...,\rho_{q})$.  Then, the special integral $\theta _{k} $ of the multidimensional Terry's problem diverges for the natural $k\ge q $ such, that
\[2kr\le v+\sum _{j=1}^{N}\sum _{i=1}^{r}k_{ji} .\]

\textbf{Theorem 3. }Let the polynomial (2) does't be represented as a sum of two polynomials of smaller number of variables, its senior form contains all of $r$ independent variables, and the matrix of exponents (3) has a rank $\rho $, $1\le \rho \le r$, is a matrix of the structure $(\rho_{1},...,\rho_{q})$. Then, the special integral of multidimensional Terry's problem   converges for all natural $k\ge q$  such that $2kr \ge 2N+r$ and
\[2kr>r+\sum _{j=1}^{N}\sum _{i=1}^{r}k_{ji}   .\]

Taking the most weak case of $\rho=1$ we get following condition of convergence doesn't dependent from singularities.

\textit{Consequence.} Let the polynomial is defined by (2), can't be represented as a sum of polynomials of smaller number of variables with common components, and their senior form contains all of independent variables. Then the special integral $\theta_{k}$ converges when $k\ge N$ and
\[2kr>r+\sum _{j=1}^{N}\sum _{i=1}^{r}k_{ji}   .\]

The consequence gives almost exact value for the convergence exponent when $2kN \le r+\sum _{j=1}^{N}\sum _{i=1}^{r}k_{ji}$.

\textbf{Theorem 4. }Let the polynomial (2) doesn't be represented as a sum of two polynomials of smaller number of variables, its senior form contains all of $r$ independent variables, and the matrix of exponents (3) has a rank $\rho $, $1\le \rho \le r$, and is a matrix of the structure $(\rho,\rho,..., \rho)$. If, in addition, the conditions of the consequence of the lemma 1 are satisfied, then the special integral of multidimensional Terry's problem  converges for all natural $k$  such that $k\rho\geq N, 2kr \ge 2N+r$ and
\[2kr>r+\sum _{j=1}^{N}\sum _{i=1}^{r}k_{ji}   .\]

The proof of the theorems 2, 3 and 4 can be spent as well as theorems 2 and 3 of the work [14], using the result of the lemma 3.

\begin{center}
\noindent \textbf{References}
\end{center}

 1.  Arnold V. I., Varchenko A. N., Huseyn-zade S. M. Singularities of differentiable mappings. M., Nauka, 2009.\textbf{}

2. Landau E. Introduction to Differential and Integral calculus. -2 ed. -М: KomKniga, 2005.

3. Arkhipov G. I., Karatsuba A. A., Chubarikov V. N. «Trigonometric integrals», Izv. Academy of Sciences. of  USSR, math. ser. (1979), v.43, №5, pp.971-1003.

4. Arkhipov G. I., Karatsuba A. A., Chubarikov V. N. «Theory of multiple trigonometric sums». M., Nauka, 1987.

5. Chubarikov V. N. «On multiple trigonometric integrals». Dokl. Academy of Sciences of USSR, (1976), v.227, №6, pp.1308-1310.

6. Chubarikov V. N. «On multiple rational trigonometric sums and multiple integrals», Mat. Notes, v.20, №1, (1976), pp. 61-68.

7. Hua Loo Keng «On the number of solutions of Tarry's problem» Acta Sci. Sinica, (1952), v.1, №1, pp. 1-76.

 8. Arkhipov G. I., Karatsuba A. A., Chubarikov V. N. «Multiple trigonometric sums and their applications», Proc. MIAN, (1980), v.151, pp.1-128.\textbf{}

9. Jabbarov I. Sh. On the structure of some algebraic varieties. Transactions of NAS of Azerbaijan, Issue Mathematics, 36 (1), 74-82 (2016). Series of Physical-Technical and Mathematical Sciences.

10. Dzhabbarov I. Sh. ``On an identity of Harmonic Analysis and its applications''. Proceedings of AS USSR, 1990, v, 314, №5, p. 1052-1054.

11. Dzhabbarov I. Sh. ``On estimates of trigonometrically integrals» Transactions of RAS, 1994, v.207, pp. 82-92.

12. Dzhabbarov I. Sh. ``On estimates of trigonometrically integrals» Chebishevskii sbornik,  v. 11, issue 1(2010), pp. 85-108.

13. Dzhabbarov I. Sh. ``On convergence exponent of the special integral of two dimensional Terry's problem'' Scientific Notes of Orlov State University, №6 (50), 2012, pp. 80-89.

14. Dzhabbarov I. Sh. ``On convergence exponent of the special integral of multi-dimensional Terry's problem'' Chebishevskii sbornik,  v.14, issue 2(2013), pp. 74-103.

15. Jabbarov I. Sh. “On convergence exponent of the special integral of two dimensional Terry’s problem”. arXiv:1302.3888v2{math.NT] 11 Mar 2017.

16. Shilov G. E Mathematical analysis. Functions of several real variables M., Nauka, 1972.

17. I. A. Ikromov. On the convergence exponent of trigonometric integrals, Tr. Mat. Inst.
Steklova, 1997, Volume 218, 179–189.

18. M. A. Chahkiev. On the convergence exponent of the singular integral in the multidimesio
nal analogue of Tarry’s problem, Izv. RAN. Ser. Mat., 2003, Volume 67, Issue 2, 211–224.
	
19. W. F. Trench. Introduction to real analysis. Free Edition 1.05, November 2010., 574p.

\noindent

\noindent

\end{document}